\documentclass[amstex,12pt,thmsa,sw20bams]{article}%
\usepackage{amsfonts}
\usepackage[utf8]{inputenc}
\usepackage{graphicx}
\usepackage[francais]{babel}
\usepackage{amsmath}
\usepackage{amssymb}%
\providecommand{\U}[1]{\protect\rule{.1in}{.1in}}
\ifx\pdfoutput\relax\let\pdfoutput=\undefined\fi
\newcount\msipdfoutput
\ifx\pdfoutput\undefined\else
\ifcase\pdfoutput\else
\ifx\paperwidth\undefined\else
\ifdim\paperheight=0pt\relax\else\pdfpageheight\paperheight\fi
\ifdim\paperwidth=0pt\relax\else\pdfpagewidth\paperwidth\fi
\fi\fi\fi
\begin{document}

\author{Choulakian V, Universit\'{e} de Moncton, Canada
\and vartan.choulakian@umoncton.ca}
\title{Geometric Interpretation of the Robustness of\\Taxicab Singular Value Decomposition of\\High-Dimensional Sparse Data}
\date{July 2026}
\maketitle

\begin{abstract}
This paper examines the fundamental properties of projections in L2{}-normed
(Euclidean) and L1-normed (Taxicab) spaces through two complementary
frameworks: projection operators and distance minimization. The problem may be
viewed as the most elementary regression problem, namely the projection of one
point onto another. Its relevance stems from the fact that any SVD-like matrix
decomposition can be interpreted as a pair of simultaneous regressions of the
rows and columns of a data matrix within either Euclidean or Taxicab geometry.
Whereas orthogonality (L2-conjugacy) is the key concept in Euclidean geometry,
L1-conjugacy, characterized by the sign function, assumes an analogous role in
Taxicab geometry. This contrast is further illustrated by the relationship
between the classical Pythagorean theorem and its Taxicab counterpart. The
primary goal of this study is to compare three decompositions: the classical
singular value decomposition (SVD), the Taxicab singular value decomposition
(TSVD), and the L1-min SVD. We also provide a geometric explanation for the
robustness of TSVD in the analysis of extremely sparse, high-dimensional
datasets by comparing the volumes of hyperspheres inscribed in hypercubes
under the two geometries. .



\end{abstract}

\section{\textbf{Introduction 1}}

This paper compares two fundamental types of projections in $l_{p}$- normed
spaces for $p=1$ and $2$: operator projections and metric projections based on
distance minimization. We emphasize that both notions of projection in Banach
spaces remain active areas of research in functional analysis; see, for
example, Agniel (2021) [$\mathbf{1}$] and Alber (1996, 2017) [$\mathbf{2,3}$].
It is well known that in Euclidean geometry, i.e., in $l_{2}$- normed spaces,
operator and metric projections coincide, whereas in Taxicab geometry, i.e.,
in $l_{1}$-normed spaces, they differ. This divergence provides the core
mathematical explanation for why $l_{1}$-based metric projection methods are
robust to outliers, while $l_{1}$-based projection operators behave
differently when applied to high-dimensional sparse data.

We note that $l_{1}$- normed space is named Taxicab geometry by Karl Menger in
1952, see among others Golland (1990) [$\mathbf{4}$].

We begin by examining the simplest mathematical--statistical problem: the
regression of one point onto another. This basic setting is fundamental, as it
provides a framework for interpreting the Singular Value Decomposition (SVD)
as a pair of simultaneous regressions involving the rows and columns of a
matrix. Our main objective is to compare three matrix decompositions: the
classical Euclidean SVD of Eckart and Young (1936) [$\mathbf{5}$], the Taxicab
SVD (TSVD) introduced by Choulakian (2006) [$\mathbf{6}$], and the $l_{1}%
\min-$ SVD proposed by Hawkins et al. (2001) [$\mathbf{7}$]. Through this
comparison, we offer a new geometric perspective that helps explain the
empirical robustness of TSVD in high-dimensional sparse data visualization;
see, for example, Choulakian (2017) [$\mathbf{8}$], Choulakian and
Allard(2025) [$\mathbf{9}$]), Gauthier and Choulakian (2015) [$\mathbf{10}$],
Choulakian and Tibeiro (2013) [$\mathbf{11}$], Choulakian et al. (2006, 2013,
2023) [$\mathbf{12,13,14}$]. Additionally, we attempt to show that for
extremely sparse datasets, both ordinary and $l_{1}\min-$ SVD behave in a
nearly identical manner. Specifically, we examine how the volume of an
inscribed hypersphere scales relative to that of its enclosing hypercube in
high dimensions; see, among others, Scott (2015, section 1.5) [$\mathbf{15}$]
and Zaki and Meira\ (2014, chapter 6) [$\mathbf{16}$], who cite Kendall (1961)
[$\mathbf{17}$] as primary reference.

The paper is structured as follows: Section 2 reviews $l_{p}$- normed spaces
for $p=1,2$. Sections 3 and 4 contrast the Euclidean and Taxicab variations of
the Pythagorean Theorem, and section 5 presents a numerical example. Section 6
introduces the three SVD variants and presents a numerical example of an
extremely sparse dataset; section 7 develops a high-dimensional geometric
interpretation of TSVD robustness. Finally, section 8 concludes the paper.

\section{Preliminaries on $l_{p}$- normed spaces for $p=1$ and 2}

The $l_{p}$ - norm of a vector $\mathbf{x}=(x_{i})\in R^{n}$ is $||\mathbf{x}%
||_{p}=(\sum_{i=1}^{n}|x_{i}|_{p})^{1/p}$ for $p\geq1$ and $||\mathbf{x}%
||_{\infty}=\max_{i}|x_{i}|.$ We note: $l_{p}^{n}:=(R^{n},||.||_{p})$ the
finite dimensional Banach space (named also Minkowski space); that is, $R^{n}$
is the $n-$dimensional complete vector space with the $l_{p}$ - norm,
$||.||_{p}$, for $p\geq1$.

A pair of numbers $(p,q)$ are conjugate if $1/p+1/q=1$ for $p\geq1$ and
$q\geq1$; this means that the normed space $l_{p}^{n}$ and its dual $l_{q}%
^{n}$ are related; in particular $||\mathbf{x}||_{p}=\varphi(\mathbf{x}%
)^{\prime}$ $\mathbf{x}$, where the norming functional $\varphi(\mathbf{x})\in
l_{q}^{n}$ and $||\varphi(\mathbf{x})||_{q}=1;$ for further details see
Choulakian (2016) [$\mathbf{18}$]. In this paper we focus on $(p,q)=(2,2)$ and
$(1,\infty)$. We have:

a) $||\mathbf{x}||_{2}=\varphi(\mathbf{x})^{\prime}$ $\mathbf{x}$, where
$\varphi(\mathbf{x})=\mathbf{x}/||\mathbf{x}||_{2}\in l_{2}^{n}$ and
$||\varphi(\mathbf{x})||_{2}=1.$

b) $||\mathbf{x}||_{1}=\varphi(\mathbf{x})^{\prime}$ $\mathbf{x}$, where
$\varphi(\mathbf{x})=sgn(\mathbf{x)=(}sign(x_{i}))\in l_{\infty}^{n}$ and
$||\varphi(\mathbf{x})||_{\infty}=1.$

c) Projection operator onto a vector $\mathbf{x\in}l_{p}^{n}$ is the matrix
$\mathbf{Q}_{\mathbf{x}}\mathbf{=}\frac{\mathbf{x}}{||\mathbf{x}||_{p}}%
\varphi(\mathbf{x})^{\prime}$, which is idempotent
\[
\mathbf{Q}_{\mathbf{x}}\mathbf{Q}_{\mathbf{x}}\mathbf{x=Qx=x.}%
\]

The following simple result is basic.\bigskip

\textbf{Lemma 1}: a) Let $\mathbf{x}\in l_{p}^{n}$ and $\mathbf{Q}%
_{\mathbf{x}}\mathbf{=}\frac{\mathbf{x}}{||\mathbf{x}||_{p}}\varphi
(\mathbf{x})^{\prime}$; then for any $\mathbf{y}\in l_{p}^{n}$, $\varphi
\mathbf{(\mathbf{x})^{\prime}(I\mathbf{-}Q}_{\mathbf{x}}\mathbf{)y=\ }%
0\mathbf{.}$

\textit{proof}: We have:%
\begin{align*}
\varphi\mathbf{(\mathbf{x})^{\prime}(I\mathbf{-}Q}_{\mathbf{x}}\mathbf{)y}  &
=\varphi\mathbf{(\mathbf{x})^{\prime}[\mathbf{y-}\frac{\mathbf{x}%
}{||\mathbf{x}||_{p}}}\varphi\mathbf{(\mathbf{x})^{\prime}y]}\\
&  =\varphi\mathbf{(\mathbf{x})^{\prime}\mathbf{y-}\frac{||\mathbf{x}||_{p}%
}{||\mathbf{x}||_{p}}}\varphi\mathbf{(\mathbf{x})^{\prime}y)}\\
&  =0.
\end{align*}

b) $\varphi(\mathbf{x})^{\prime}\mathbf{y}=0$ is equivalent to $\varphi
(\mathbf{x})^{\prime}\mathbf{Q}_{\mathbf{y}}=\mathbf{0.}$\bigskip

We note that $\mathbf{Q}_{\mathbf{x}}\mathbf{y}\in l_{p}^{n}$ is the
projection of $\mathbf{y}$ onto $\mathbf{x;}$ $\mathbf{(I\mathbf{-}%
Q}_{\mathbf{x}}\mathbf{)y}$ is the residual vector. The central aim of this
paper is to discuss the relations between $\mathbf{y,}$ and $\mathbf{Q}%
_{\mathbf{x}}\mathbf{y}$ and $\mathbf{(I\mathbf{-}Q}_{\mathbf{x}}\mathbf{)y.}$
For instance, by applying Lemma 1 in a Gram-Schmidt like process to a set of
$1\leq k\leq n$ independent vectors $\mathbf{x}_{i}\in l_{p}^{n}$ for
$i=1,...,k$, \textit{one gets a set of }$l_{p}$\textit{ - conjugate ordered
vectors }$\mathbf{y}_{i}\in l_{p}^{n}$\textit{ for }$i=1,...,k.$ That is,
$\varphi(\mathbf{y}_{\beta})^{\prime}\mathbf{y}_{\alpha}=0$ for $k\geq
\alpha>\beta\geq1$.\bigskip

\textbf{Example 1: }Gram-Schmidt process to a set of $k=3\leq n$ independent
vectors $\mathbf{x}_{i}\in l_{p}^{n}$ for $i=1,2,3$ and $\mathbf{Q}_{i}=$
$\frac{\mathbf{y}_{i}}{||\mathbf{y}_{i}||_{p}}\varphi(\mathbf{y}_{i})^{\prime
}$ is:

(1) $\mathbf{y}_{1}=\mathbf{x}_{1}$

(2) $\mathbf{y}_{2}=(\mathbf{I}-\mathbf{Q}_{1})\mathbf{x}_{2}.$ Then:

(2a) $\varphi(\mathbf{y}_{1})^{\prime}\mathbf{y}_{2}=\varphi(\mathbf{y}%
_{1})^{\prime}(\mathbf{I}-\mathbf{Q}_{1})\mathbf{x}_{2}=0$ by Lemma 1a; that
is, the ordered vectors $\mathbf{y}_{1}$ and $\mathbf{y}_{2}$ are $l_{p}$ -
conjugate: $\varphi(\mathbf{y}_{1})^{\prime}\mathbf{y}_{2}=0.$ Thus by Lemma
1b: $\varphi(\mathbf{y}_{1})^{\prime}\mathbf{y}_{2}=\varphi(\mathbf{y}%
_{1})^{\prime}\mathbf{Q}_{2}=\mathbf{0}.$

(3) $\mathbf{y}_{3}=(\mathbf{I}-\mathbf{Q}_{2})(\mathbf{I}-\mathbf{Q}%
_{1})\mathbf{x}_{3}.$ Then:

(3a) $\varphi(\mathbf{y}_{2})^{\prime}\mathbf{y}_{3}=\varphi(\mathbf{y}%
_{2})^{\prime}(\mathbf{I}-\mathbf{Q}_{2})[(\mathbf{I}-\mathbf{Q}%
_{1})\mathbf{x}_{3}]\ \mathbf{=\ }0$ by Lemma 1a.

(3b) We have:
\begin{align*}
\varphi(\mathbf{y}_{1})^{\prime}\mathbf{y}_{3}  &  =\varphi(\mathbf{y}%
_{1})^{\prime}(\mathbf{I}-\mathbf{Q}_{2})(\mathbf{I}-\mathbf{Q}_{1}%
)\mathbf{x}_{3}\\
&  =\varphi(\mathbf{y}_{1})^{\prime}(\mathbf{I}-\mathbf{Q}_{1})\mathbf{x}%
_{3}-\varphi(\mathbf{y}_{1})^{\prime}\mathbf{Q}_{2}[(\mathbf{I}-\mathbf{Q}%
_{1})\mathbf{x}_{3}]\\
&  =0-\varphi(\mathbf{y}_{1})^{\prime}\mathbf{Q}_{2}[(\mathbf{I}%
-\mathbf{Q}_{1})\mathbf{x}_{3}]\text{ \ by Lemma 1a}\\
&  =0-0\text{ \ by (2a) or Lemma 1b;}%
\end{align*}
that is, the ordered vectors $\mathbf{y}_{1},$ $\mathbf{y}_{2}$ and
$\mathbf{y}_{3}$ are $l_{p}$ - conjugate: $\varphi(\mathbf{y}_{\beta}%
)^{\prime}\mathbf{y}_{\alpha}=0$ for $3\geq\alpha>\beta\geq1.$

We can summarize the above computations in the following way: Let
$\mathbf{Y}=[\mathbf{y}_{1}\ \mathbf{y}_{2}\ \mathbf{y}_{3}]\in R^{n\times3}$
be a set of three $l_{p}$ - conjugate ordered vectors$.$

a) For $p=1$, the matrix $\mathbf{Y}$ is named $l_{1}$ - conjugate, for%

\[
sgn(\mathbf{Y)}^{\prime}\mathbf{Y=}\left[
\begin{array}
[c]{ccc}%
||\mathbf{y}_{1}||_{1} & 0 & 0\\
sgn(\mathbf{y}_{2}\mathbf{)}^{\prime}\mathbf{y}_{1} & ||\mathbf{y}_{2}||_{1} &
0\\
sgn(\mathbf{y}_{3}\mathbf{)}^{\prime}\mathbf{y}_{1} & sgn(\mathbf{y}%
_{3}\mathbf{)}^{\prime}\mathbf{y}_{2} & ||\mathbf{y}_{3}||_{1}%
\end{array}
\right]  .
\]

b) For $p=2$, the matrix $\mathbf{Y}$ is orthogonal ($l_{2}$ - conjugate),
for
\[
\mathbf{Y}^{\prime}\mathbf{Y=Diag(}||\mathbf{y}_{1}||_{2}^{2}\ \ ||\mathbf{y}%
_{2}||_{2}^{2}\ \ ||\mathbf{y}_{3}||_{2}^{2}).
\]
\bigskip

In the next two sections, we compare the Pythagorean Theorem and its $l_{1}$ -
norm analogue based on operator projections and metric projections based on
distance minimization.

\section{Pythagorean Theorem}

Let $\mathbf{y}\in l_{2}^{n}$, $\mathbf{x}\in l_{2}^{n}$, $\mathbf{Q}%
_{\mathbf{x}}\mathbf{=}\frac{\mathbf{x}}{||\mathbf{x}||_{2}}\varphi
(\mathbf{x})^{\prime}$ where $\varphi(\mathbf{x})=\mathbf{x}/||\mathbf{x}%
||_{2}\in l_{2}^{n}$ and $\mathbf{x}$ and $\mathbf{y}$ are linearly
independent. The orthogonal (Euclidean or $l_{2}$ - conjugate) projection of
$\mathbf{y}$ onto $\mathbf{x}$ is the vector $\widehat{\mathbf{y}}%
_{l2}=\mathbf{Q}_{\mathbf{x}}\mathbf{y}=\alpha_{eucl}\mathbf{x}$\textbf{
}where\textbf{ }%
\begin{align}
\alpha_{eucl}  &  =\frac{\mathbf{x}^{\prime}\mathbf{y}}{||\mathbf{x}||_{2}%
^{2}}\tag{1}\\
&  =\arg\min_{\alpha}||\mathbf{y-}\alpha\mathbf{x}||_{2}^{2}; \tag{2}%
\end{align}
and the residual vector is $\mathbf{(y-}\widehat{\mathbf{y}}_{l2}%
)\neq\mathbf{0}$, for $\mathbf{x}$ and $\mathbf{y}$ are linearly independent.
It is well known, or by Lemma 1, that ($\mathbf{y-}\widehat{\mathbf{y}}_{l2})$
and $\widehat{\mathbf{y}}_{l2}$ are orthogonal; thus we have
\begin{equation}
||\mathbf{y}||_{2}^{2}=||\mathbf{y-}\widehat{\mathbf{y}}_{l2}||_{2}%
^{2}+||\widehat{\mathbf{y}}_{l2}||_{2}^{2}\text{ \ \ \ \ \ if and only if
}(\mathbf{y-}\widehat{\mathbf{y}}_{l2}\mathbf{)}^{\prime}\widehat{\mathbf{y}%
}_{l2}\ \mathbf{=\ }0, \tag{3}%
\end{equation}
which is the Pythagorean Theorem; that is, the Pythagorean Theorem is based on
the orthogonality of the two vectors $(\mathbf{y-}\widehat{\mathbf{y}}_{l2})$
and $\widehat{\mathbf{y}}_{l2}$.

We also have: for two linearly independent non-orthogonal vectors $\mathbf{y}$
and $\mathbf{z}$%

\begin{align}
||\mathbf{y}||_{2}^{2}  &  =||\mathbf{y-z+z}||_{2}^{2}\nonumber\\
&  \neq||\mathbf{y-z}||_{2}^{2}+||\mathbf{z}||_{2}^{2}\text{.} \tag{4}%
\end{align}

Note that the inequality (4) does not represent the standard triangle
inequality $||\mathbf{y-z+z}||_{2}\leq||\mathbf{y-z}||_{2}+||\mathbf{z}%
||_{2}.$

\textbf{Example 2}: The aim of this example is to show the difference between
equations (4) and (6). This example concerns equation (4).

Case a) Let $\mathbf{y=\ }\left(  _{2}^{3}\right)  $ and $\mathbf{z=\ }\left(
_{1}^{2}\right)  ;$ so ($\mathbf{y-}$ $\mathbf{z)=\ }\left(  _{1}^{1}\right)
.$ Thus%
\[
||\mathbf{y}||_{2}^{2}=13>||\mathbf{z}||_{2}^{2}+||\mathbf{y-z}||_{2}%
^{2}=5+2.
\]

Case b) Let $\mathbf{y=\ }\left(  _{2}^{3}\right)  $ and $\mathbf{z=\ }\left(
_{-1}^{-2}\right)  ;$ so ($\mathbf{y-}$ $\mathbf{z)=\ }\left(  _{3}%
^{5}\right)  .$ Thus%
\[
||\mathbf{y}||_{2}^{2}=13<||\mathbf{z}||_{2}^{2}+||\mathbf{y-z}||_{2}%
^{2}=5+34.
\]

\section{Taxicab Pythagorean Theorem }

In this section we study the relations between the two vectors $\mathbf{y\in
\ }l_{1}^{n}$ and $\widehat{\mathbf{y}}_{l1}\mathbf{\in\ }l_{1}^{n}$ by noting
that the four quadrants $(sgn(y_{i}),sgn\mathbf{(}\widehat{y}_{l1}(i))$ play
fundamental role. In Proposition 1, the two opposing quadrants $(-,-)$ and
$(+,+)$ are necessary in Taxicab Pythagorean Theorem represented by equation (5).

\textbf{Proposition 1}: The $l_{1}-$ norm analogues of equations (3) and (4)
are
\begin{align}
||\mathbf{y}||_{1}  &  =||\widehat{\mathbf{y}}_{l1}||_{1}+||\mathbf{y-}%
\widehat{\mathbf{y}}_{l1}||_{1}\text{ if and only if }0\leq\frac
{\widehat{y}_{l1}(i)}{y_{i}}\leq1\mathbf{.}\tag{5}\\
||\mathbf{y}||_{1}  &  <||\widehat{\mathbf{y}}_{l1}||_{1}+||\mathbf{y-}%
\widehat{\mathbf{y}}_{l1}||_{1}\text{ \ \ \ \ \ \ \ \ \ \ \ \ otherwise.}
\tag{6}%
\end{align}

Note that the inequality (6) represents the standard triangle inequality
$||\mathbf{y-z+z}||_{1}\leq||\mathbf{y-z}||_{1}+||\mathbf{z}||_{1}.$

\textit{Proof}: Equation (5) is
\[
\sum_{i}|y_{i}|\ =\sum_{i}|y_{i}\mathbf{-}\widehat{y}_{l1}(i)|\ +\sum
_{i}|\widehat{y}_{l1}(i)|\ \text{if and only if }0\leq\frac{\widehat{y}%
_{l1}(i)}{y_{i}}\leq1\text{.}%
\]
We have for an index $i=1,...,n$:
\begin{align}
|y_{i}|\  &  <(1-\frac{\widehat{y}_{l1}(i)}{y_{i}})\ |y_{i}|\text{ \ \ \ \ if
}\frac{\widehat{y}_{l1}(i)}{y_{i}}<0\tag{7}\\
&  <|y_{i}\mathbf{-}\widehat{y}_{l1}(i)|\ +|\widehat{y}_{l1}(i)|;\nonumber
\end{align}
and%

\begin{align}
|y_{i}|\  &  <\frac{\widehat{y}_{l1}(i)}{y_{i}}\ |y_{i}|\text{ \ \ \ \ \ if
}\frac{\widehat{y}_{l1}(i)}{y_{i}}>1\tag{8}\\
&  <|y_{i}\mathbf{-}\widehat{y}_{l1}(i)|\ +|\widehat{y}_{l1}(i)|;\nonumber
\end{align}
finally based on (7) and (8)%

\begin{align}
|y_{i}|\  &  =(1-\frac{\widehat{y}_{l1}(i)}{y_{i}})\ |y_{i}|+\frac
{\widehat{y}_{l1}(i)}{y_{i}}\ |y_{i}|\text{\ \ if and only if }0\leq
\frac{\widehat{y}_{l1}(i)}{y_{i}}\leq1\tag{9}\\
&  =|y_{i}\mathbf{-}\widehat{y}_{l1}(i)|\ +|\widehat{y}_{l1}(i)|.\nonumber
\end{align}

Taxicab Pythagorean Theorem, equation (5), is based on (9); and the triangle
inequality in (6) is based on (7) or (8).\bigskip

\textbf{Remarks}

a) Let $[\mathbf{y}\ \ \widehat{\mathbf{y}}_{l1}]\in R^{n\times2};$ then
Taxicab Pythagorean Theorem, equation (5), can be interpreted in two different
ways: first, the two column vectors $\mathbf{y}\ $and$\ \widehat{\mathbf{y}%
}_{l1}$ belong to the same orthant in $R^{n}$ and $|y_{i}|\ \geq
|\widehat{y}_{l1}(i)|$; second, each row vector $(y_{i},\widehat{y}%
_{l1}(i))\in quadrant(+,+)$ or $quadrant(-,-)$ in $R^{2}$ and $|y_{i}%
|\ \geq|\widehat{y}_{l1}(i)|.$

b) Note that equations (1) and (2) are related with equation (3). In the next
two subsections, we study the corresponding relations in the Taxicab geometry.

\subsection{Regression by l$_{1}-$norm projection operator\ }

Based on the $l_{1}$ - norm variant of equation (1) we get:\bigskip

\textbf{Corollary 1}: Let $\mathbf{y}\in l_{1}^{n}$, $\mathbf{x}\in l_{1}%
^{n},$ $\mathbf{Q}_{\mathbf{x}}\mathbf{=}\frac{\mathbf{x}}{||\mathbf{x}||_{1}%
}sgn(\mathbf{x)}^{\prime}$ and $\mathbf{x}$\textbf{ }and $\mathbf{y}$\textbf{
}are linearly independent. The $l_{1}$ - norm (Taxicab) projection of
$\mathbf{y}$ onto $\mathbf{x}$ is the vector $\widehat{\mathbf{y}}%
_{l1proj}=\alpha_{l1proj}\mathbf{x}$, where\textbf{ }%
\begin{equation}
\alpha_{l1proj}=\frac{sgn(\mathbf{x)}^{\prime}\mathbf{y}}{||\mathbf{x}||_{1}}
\tag{10}%
\end{equation}
and the residual vector is $\mathbf{y-}\widehat{\mathbf{y}}_{l1proj}%
\neq\mathbf{0}.$ Then inequality (6) is true%
\begin{equation}
||\mathbf{y}||_{1}<||\widehat{\mathbf{y}}_{l1proj}||_{1}+||\mathbf{y-}%
\widehat{\mathbf{y}}_{l1proj}||_{1}\text{.} \tag{6}%
\end{equation}

\textit{Proof}: We have%
\begin{align}
||\widehat{\mathbf{y}}_{l1proj}||_{1}  &  =||\mathbf{Q}_{\mathbf{x}}%
\mathbf{y}||_{1}=||\alpha_{l1proj}\mathbf{x||}_{1}\nonumber\\
&  =|sgn(\mathbf{x)}^{\prime}\mathbf{y|}\text{ \ by (10)}\mathbf{.} \tag{11}%
\end{align}
We distinguish two cases:

Case a) $|sgn(\mathbf{x)}^{\prime}sgn(\mathbf{y)|\ =\ }|sgn(\mathbf{y)}%
^{\prime}sgn(\mathbf{y)|\ =\ }n$;$\ $ so
\begin{align*}
||\mathbf{y}||_{1}  &  =|sgn(\mathbf{y)}^{\prime}\mathbf{y|}\\
&  \mathbf{=|}sgn\mathbf{(x\mathbf{)}^{\prime}\mathbf{y|\ =}\ }%
||\widehat{\mathbf{y}}_{l1proj}||_{1}\text{ \ by (11)\textbf{;}}%
\end{align*}
given that $\mathbf{y}-\widehat{\mathbf{y}}_{l1proj}\ \mathbf{=y}%
-\alpha_{l1proj}\mathbf{x\neq0}$ for $\mathbf{x}$\textbf{ }and $\mathbf{y}%
$\textbf{ }are linearly independent, thus inequality (6) is true.

Case b) $|sgn(\mathbf{x)}^{\prime}sgn(\mathbf{y)|\ <\ }n$; which implies the
existence of two terms $sgn(x_{i_{1}}\mathbf{)\ }sgn\ (y_{i_{1}}%
\mathbf{)=-\ }sgn(x_{i_{2}}\mathbf{)\ }sgn\ (y_{i_{2}}\mathbf{)=\ }1.$

If $\alpha_{l1proj}<0,$ then $\frac{\widehat{y}_{l1}(i_{1})=\alpha
_{l1proj}\ x_{i_{1}}}{y_{i_{1}}}<0,$ which implies inequality (6) is true.

If $\alpha_{l1proj}>0,$ then $\frac{\widehat{y}_{l1}(i_{2})=\alpha
_{l1proj}\ x_{i_{2}}}{y_{i_{2}}}<0,$ which implies inequality (6) is true.

We note that in this case,
\begin{align*}
||\mathbf{y}||_{1}  &  =|sgn(\mathbf{y)}^{\prime}\mathbf{y|}\\
&  \mathbf{>|}sgn\mathbf{(x\mathbf{)}^{\prime}\mathbf{y|}\ =}\text{
}||\widehat{\mathbf{y}}_{l1proj}||_{1}\text{ \ \ \ \ \ \ \ \ \ \ by
(11)\textbf{.}}%
\end{align*}

\textbf{Remark}: $\alpha_{eucl}$ in equation (1) and $\alpha_{l1proj}$ in
equation (10) can be represented as%
\[
\alpha_{lp-proj}=\varphi(\mathbf{x})^{\prime}\mathbf{y/||x||}_{p}\text{ \ for
}lp=l_{1}\text{ and }l_{2}\text{.}%
\]
\bigskip

\subsection{Regression by l$_{1\min}-\ $distance projection }

The $l_{1}$ - norm variant of equation (2) is
\begin{align}
\alpha_{l1\min}  &  =\arg\min_{\alpha}||\mathbf{y-}\alpha\mathbf{x}%
||_{1}\tag{12}\\
&  =\arg\min_{\alpha}\sum_{i}\ |x_{i}|\ |y_{i}/x_{i}\mathbf{-}\alpha
|;\nonumber
\end{align}
the estimate $\alpha_{l1\min}$ is the weighted median of the vector
$(y_{i}/x_{i})$ with weight vector $(|x_{i}|).$ With great probability, it
follows that there exists an index $i=imin$ such that
\begin{equation}
|\frac{y_{i}}{x_{i}}\mathbf{-}\alpha_{l1\min}|\ \geq|\frac{y_{i\min}}%
{x_{i\min}}\mathbf{-}\alpha_{l1\min}|\ =0\text{ for }imin\neq i\text{ for
}i=1,...,n\text{.} \tag{13}%
\end{equation}

\textbf{Corollary 2}: Let $\mathbf{y}\in l_{1}^{n}$, $\mathbf{x}\in l_{1}%
^{n},$ and $\mathbf{x}$\textbf{ }and $\mathbf{y}$\textbf{ }are linearly
independent. The $l_{1\min}$ - norm (Taxicab) metric projection of
$\mathbf{y}$ onto $\mathbf{x}$ by (12) is the vector $\widehat{\mathbf{y}%
}_{l1\min}=\alpha_{l1\min}\mathbf{x}$, and the residual vector is
$\mathbf{y-}\widehat{\mathbf{y}}_{l1\min}\neq\mathbf{0}.$ Then by Proposition
1%
\begin{align}
||\mathbf{y}||_{1}  &  =||\widehat{\mathbf{y}}_{l1\min}||_{1}+||\mathbf{y-}%
\widehat{\mathbf{y}}_{l1\min}||_{1}\text{ if and only if }0\leq\frac
{\widehat{y}_{l1\min}(i)}{y_{i}}\leq1\mathbf{.}\tag{5}\\
||\mathbf{y}||_{1}  &  <||\widehat{\mathbf{y}}_{l1}||_{1}+||\mathbf{y-}%
\widehat{\mathbf{y}}_{l1}||_{1}\text{ \ \ \ \ \ \ \ \ \ \ \ \ otherwise.}
\tag{6}%
\end{align}

\section{\textbf{Example 3 }}

The famous recreational mathematician Gardner (1997, ch. 10)[$\mathbf{19}$]
provides many examples of Taxicab plane geometry, which distinguishes it from
the Euclidean plane geometry, and shows its flexibility; a simple example
provided by Gardner shows the difference between the Pythagorean Theorem (3)
and its Taxicab analogue (5):\textit{ "A taxicab scalene triangle with corners
A, B, and C and sides of 14, 8, and 6 is shown at the left in Figure 65. The
sides of taxi polygons must of course be taxi paths, and the paths that make
up a polygon of specified dimensions may vary in shape but not in length.
Observe how the triangle in the illustration violates the Euclidean theorem
that the sum of any two sides of a triangle must be greater than the third
side. In this case the sum of two sides equals the third: 6 + 8 equals 14."}

\subsection{Solution by Taxicab Pythagorean Theorem}

A simple generalization of Gardner's example is: Let $\mathbf{a=\ }\left(
_{y}^{x}\right)  ,$ $\mathbf{c=\ }\left(  _{\beta y}^{\alpha x}\right)  $ and
$\mathbf{a-c=\ }\left(  _{y-\beta y}^{x-\alpha x}\right)  .$ Furthermore,
$||\mathbf{a}||_{p}=14,\ ||\mathbf{c}||_{p}=6,\ ||\mathbf{a-c}||_{p}=8$ for
$p=1,2$; that is%

\begin{equation}
||\mathbf{a}||_{p}=||\mathbf{c}||_{p}+||\mathbf{a-c}||_{p}. \tag{14}%
\end{equation}

For $p=2$, (14) is true if $\mathbf{c=a(}6/14\mathbf{)}$; that is,
$\mathbf{c}$ and $\mathbf{a}$ are linearly dependent, which contradicts the
assumption of existence of a triangle.

For $p=1$, (14) is equivalent to a system of three equations in four
unknowns,
\begin{align*}
14  &  =|x|+|y|\ =||\mathbf{a}||_{1}\\
6  &  =|\alpha|\ |x|+|\beta|\ |y|\ =||\mathbf{c}||_{1}\\
8  &  =|1-\alpha|\ |x|\ +|1-\beta|\ |y|)=||\mathbf{a-c}||_{1};
\end{align*}
which has infinite number of solutions. The solutions are found in the
rectangle constructed by the four corners $\mathbf{a=\ }\left(  _{y>0}%
^{x>0}\right)  ,$ $\left(  _{y=0}^{x=0}\right)  ,\ \left(  _{y=0}%
^{x>0}\right)  $ and $\left(  _{y>0}^{x=0}\right)  $. In Gardner's problem:
$\mathbf{a=\ }\left(  _{y=8}^{x=6}\right)  $, $\mathbf{c=\ }\left(  _{\beta
y=2}^{\alpha x=4}\right)  $ and $\mathbf{a-c=\ }\left(  _{y-\beta
y=6}^{x-\alpha x=2}\right)  $; so $0\leq\alpha x/x=\alpha=4/6\leq1$ and
$0\leq\beta=2/8\leq1$, which represent the conditions, $0\leq\frac
{\widehat{y}_{l1}(i)}{y_{i}}\leq1,$ for the the validity of Taxicab
Pythagorean Theorem, equation (5).

Let us compute the projections of $\mathbf{a}$ onto $\mathbf{c}$ by the three
methods discussed above, where $\widehat{\mathbf{a}}=\theta\mathbf{c}$.

First, by Euclidean projection, equation (1): $\theta_{eucl}=\frac
{\mathbf{c}^{\prime}\mathbf{a}}{||\mathbf{c}||_{2}^{2}}=\frac{24+16}%
{16+4}=\frac{40}{20}=2.$ Thus%
\begin{align*}
||\mathbf{a}||_{2}^{2}  &  =||\theta_{eucl}\mathbf{c}||_{2}^{2}+||\mathbf{a}%
-\theta_{eucl}\mathbf{c}||_{2}^{2}\\
||\left(  _{8}^{6}\right)  ||_{2}^{2}  &  =||2\left(  _{2}^{4}\right)
||_{2}^{2}+||\left(  _{8}^{6}\right)  -2\left(  _{2}^{4}\right)  ||_{2}^{2}\\
100  &  =80+20.
\end{align*}

Second, by Taxicab projection operator, equation (10): $\theta_{l1proj}%
=\frac{sgn(c\mathbf{)}^{\prime}\mathbf{a}}{||\mathbf{c}||_{1}}=\frac{6+8}%
{4+2}=\frac{14}{6}.$ Thus%
\begin{align*}
||\mathbf{a}||_{1}  &  <||\theta_{l1proj}\mathbf{c}||_{1}+||\mathbf{a}%
-\theta_{l1proj}\mathbf{c}||_{1}\\
||\left(  _{8}^{6}\right)  ||_{1}  &  <\frac{14}{6}||\left(  _{2}^{4}\right)
||_{1}+||\left(  _{8}^{6}\right)  -\frac{14}{6}\left(  _{2}^{4}\right)
||_{1}\\
14  &  <14+|6-28/3|+|8-14/3|.
\end{align*}

Third, by $l_{1\min}$-distance minimization, equation (12):
\begin{align*}
\theta_{l1\min}  &  =\arg\min_{\theta}\sum_{i}\ |c_{i}|\ |a_{i}/c_{i}%
\mathbf{-}\theta|\\
&  =\arg\min_{\theta}[4|1.5-\theta|+2|4-\theta|]\\
&  =1.5.
\end{align*}
Thus%
\begin{align*}
||\mathbf{a}||_{1}  &  =||\theta_{l1\min}\mathbf{c}||_{1}+||\mathbf{a}%
-\theta_{l1\min}\mathbf{c}||_{1}\\
||\left(  _{8}^{6}\right)  ||_{1}  &  =1.5||\left(  _{2}^{4}\right)
||_{1}+||\left(  _{8}^{6}\right)  -1.5\left(  _{2}^{4}\right)  ||_{1}\\
14  &  =9+0+5.
\end{align*}

\subsection{Modification }

Let us modify Gardner's problem and consider: $\mathbf{a=\ }\left(
_{y=8}^{x=6}\right)  $ and $\mathbf{c=\ }\left(  _{\beta y=-2}^{\alpha
x=4}\right)  ,$ from which it follows $\mathbf{a-c=\ }\left(  _{y-\beta
y=10}^{x-\alpha x=2}\right)  .$ We\ get the inequality (6), because $\beta
y/y=\beta=-2/8<0$:
\begin{align}
||\mathbf{a}||_{1}  &  <||\mathbf{c}||_{1}+||\mathbf{a-c}||_{1}\tag{6}\\
14  &  <6+12.\nonumber
\end{align}
\ 

Let us compute the projections of $\mathbf{a}$ on $\mathbf{c}$ by the three
methods discussed above, where $\widehat{\mathbf{a}}=\theta\mathbf{c}$.

First, Euclidean projection by (1): $\theta_{eucl}=\frac{\mathbf{c}^{\prime
}\mathbf{a}}{||\mathbf{c}||_{2}^{2}}=\frac{24-16}{16+4}=\frac{8}{20}.$

Second, Taxicab projection operator by (10): $\theta_{l1proj}=\frac
{sgn(c\mathbf{)}^{\prime}\mathbf{a}}{||\mathbf{c}||_{1}}=\frac{6-8}{4+2}%
=\frac{-2}{6}.$ Thus%
\begin{align*}
||\mathbf{a}||_{1}  &  <||\theta_{l1proj}\mathbf{c}||_{1}+||\mathbf{a}%
-\theta_{l1proj}\mathbf{c}||_{1}\\
||\left(  _{8}^{6}\right)  ||_{1}  &  <\frac{1}{3}||\left(  _{-2}^{4}\right)
||_{1}+||\left(  _{8}^{6}\right)  +\frac{1}{3}\left(  _{-2}^{4}\right)
||_{1}\\
14  &  <2+|6+4/3|+|8-2/3|.
\end{align*}

Third, by $l_{1}$-distance minimization by (12):
\begin{align*}
\theta_{l1\min}  &  =\arg\min_{\theta}\sum_{i}\ |c_{i}|\ |a_{i}/c_{i}%
\mathbf{-}\theta|\\
&  =\arg\min_{\theta}\ [2|-4-\theta|+4|1.5-\theta|]\\
&  =1.5.
\end{align*}

Note that $\theta_{l1\min}=1.5$ in both cases; however Taxicab Pythagorean
Theorem (5) is not satisfied in this case:
\begin{align}
||\mathbf{a}||_{1}  &  <||\theta_{l1\min}\mathbf{c}||_{1}+||\mathbf{a}%
-\theta_{l1\min}\mathbf{c}||_{1}\tag{6}\\
||\left(  _{8}^{6}\right)  ||_{1}  &  <1.5||\left(  _{-2}^{4}\right)
||_{1}+||\left(  _{8}^{6}\right)  -1.5\left(  _{-2}^{4}\right)  ||_{1}%
\nonumber\\
14  &  <9+(0+11).\nonumber
\end{align}

\section{SVD, Taxicab SVD and l$_{1}\min-$ SVD }

Eckart and Young (1936)[$\mathbf{5}$] showed that Pearson (1901)
[$\mathbf{20}$] and Hotelling (1933) [$\mathbf{21}$] approaches of principal
component analysis (PCA) of a dataset \textbf{X} is based on singular value
decomposition (SVD); which is the workhorse for classical multidimensional
data analysis based on Euclidean geometry, see Jolliffe (2002) [$\mathbf{22}%
$]. Similarly, Choulakian (2003) [$\mathbf{23}$] showed that Taxicab SVD
(TSVD) of a variance-covariance matrix \textbf{S} is equivalent to Burt
(1917)[$\mathbf{24}$]-Thurstone (1931)[$\mathbf{25}$] approach of centroid
decomposition of the covariance matrix \textbf{S}. We attempt to see the
similarities and differences among the three SVDs based on the three
projections discussed above.

Let $\mathbf{X=(}x_{ij})$\ be a matrix of size $I\times J$ and
$rank(\mathbf{X})=k$\textbf{.} The three SVDs are of the form
\begin{equation}
\mathbf{X}=\sum_{\alpha=1}^{k}\mathbf{a}_{\alpha}\mathbf{b}_{\alpha}^{\prime
}/\delta_{\alpha}, \tag{15}%
\end{equation}
where
\begin{align*}
\delta_{\alpha}  &  =||\mathbf{a}_{\alpha}||_{p}=\varphi(\mathbf{a}_{\alpha
})^{\prime}\mathbf{a}_{\alpha}\\
&  \mathbf{=\ }||\mathbf{b}_{\alpha}||_{p}=\varphi(\mathbf{b}_{\alpha
})^{\prime}\mathbf{b}_{\alpha}.
\end{align*}
The parameters $(\delta_{\alpha},\mathbf{a}_{\alpha},\mathbf{b}_{\alpha})$ for
$\alpha=1,...,k$ are generally computed in sequential, stepwise manner. Let
$\mathbf{X}_{1}=\mathbf{X=(}x_{ij})$ and
\begin{align*}
\mathbf{X}_{\alpha+1}  &  =\mathbf{X}_{1}-\sum_{\beta=1}^{\alpha}%
\mathbf{a}_{\beta}\mathbf{b}_{\beta}^{\prime}/\delta_{\beta}\text{
\ for\ \ }\alpha=1,...,k-1\\
&  =\mathbf{X}_{\alpha}-\mathbf{a}_{\alpha}\mathbf{b}_{\alpha}^{\prime}%
/\delta_{\alpha}%
\end{align*}
be the residual matrix at the $\alpha$th iteration. Let $||\mathbf{X}%
_{\substack{\alpha\\}}||_{p}^{p}=\sum_{i,j}|x_{\substack{\alpha\\}%
}(i,j)|_{\text{ }}^{p}$ for $p=1$ and $2$; then by triangle inequality
\begin{align*}
||\mathbf{X}_{\alpha}||_{p}  &  =||\mathbf{a}_{\alpha}\mathbf{b}_{\alpha
}^{\prime}/\delta_{\alpha}+\mathbf{X}_{\alpha}-\mathbf{a}_{\alpha}%
\mathbf{b}_{\alpha}^{\prime}/\delta_{\alpha}||_{p}\\
&  \leq||\mathbf{a}_{\alpha}\mathbf{b}_{\alpha}^{\prime}||_{p}/\delta_{\alpha
}+||\mathbf{X}_{\substack{\alpha\\}}-\mathbf{a}_{\alpha}\mathbf{b}_{\alpha
}^{\prime}/\delta_{\alpha}||_{p}\text{ for }p=1\ \text{and}\ 2.
\end{align*}
The term $\mathbf{a}_{\alpha}\mathbf{b}_{\alpha}^{\prime}/\delta_{\alpha}$ can
be calculated either by maximizing $||\mathbf{a}_{\alpha}\mathbf{b}_{\alpha
}^{\prime}||_{p}^{p}=\sum_{i}\sum_{j}|a_{i}b_{j}|^{p}=\delta_{\alpha}^{2p}$ or
by minimizing $||\mathbf{X}_{\alpha}-\mathbf{a}_{\alpha}\mathbf{b}_{\alpha
}^{\prime}/\delta_{\alpha}||_{p}^{p}$ for $p=1$ or $2$, that we discuss in the
next two subsections.

\subsection{SVD and Taxicab SVD \ }

Hotelling (1933)[$\mathbf{21}$] named principal component analysis (PCA) the
decomposition (15) by maximizing the spectral (Euclidean) norm; while Burt
(1917)[$\mathbf{24}$]-Thurstone (1931)[$\mathbf{25}$] named centroid
decomposition the decomposition (15) of a covariance matrix $\mathbf{X}%
^{\prime}\mathbf{X}$ by maximizing the Grothendieck norm, see Choulakian
(2006, 2016)[$\mathbf{6,18}$]; from which the $l_{p}$-norm projection
operators are derived.

The spectral (Euclidean) norm of a matrix $\mathbf{X}_{\alpha}$ is%
\begin{align}
||\mathbf{X}_{\alpha}||_{2\rightarrow2}  &  =\max_{\mathbf{b}_{\alpha}}%
\frac{||\mathbf{X}_{\alpha}\mathbf{b}_{\alpha}||_{2}}{||\mathbf{b}_{\alpha
}||_{2}}\tag{16}\\
&  =\max_{\mathbf{a}_{\alpha}}\frac{||\mathbf{X}_{\alpha}^{\prime}%
\mathbf{a}_{\alpha}||_{2}}{||\mathbf{a}_{\alpha}||_{2}}\nonumber\\
&  =\max_{\mathbf{a}_{\alpha},\mathbf{b}_{\alpha}}\mathbf{a}_{\alpha}^{\prime
}\mathbf{X}_{\alpha}\mathbf{b}_{\alpha}=\delta_{\alpha}\text{ subject to
}||\mathbf{a}_{\alpha}||_{2}=||\mathbf{b}_{\alpha}||_{2}=1.\nonumber
\end{align}

The Grothendieck norm of a matrix $\mathbf{X}_{\alpha}$ is%
\begin{align}
||\mathbf{X}_{\alpha}||_{\infty\rightarrow1}  &  =\max_{\varphi(\mathbf{b}%
_{\alpha})}\frac{||\mathbf{X}_{\alpha}\varphi(\mathbf{b}_{\alpha})||_{1}%
}{||\varphi(\mathbf{b}_{\alpha})||_{\infty}}=\max_{\varphi(\mathbf{b}_{\alpha
})\in\left\{  -1,1\right\}  ^{J}}||\mathbf{X}_{\alpha}\varphi(\mathbf{b}%
_{\alpha})||_{1}\nonumber\\
&  =\max_{\varphi(\mathbf{a}_{\alpha})}\frac{||\mathbf{X}_{\alpha}^{\prime
}\varphi(\mathbf{a}_{\alpha})||_{1}}{||\varphi(\mathbf{a}_{\alpha})||_{\infty
}}=\max_{\varphi(\mathbf{a}_{\alpha})\in\left\{  -1,1\right\}  ^{I}%
}||\mathbf{X}_{\alpha}^{\prime}\varphi(\mathbf{a}_{\alpha})||_{1}\nonumber\\
&  =\max_{||\varphi(\mathbf{a}_{\alpha})||_{\infty}=||\varphi(\mathbf{b}%
_{\alpha})||_{\infty}=1}\varphi(\mathbf{a}_{\alpha})^{\prime}\mathbf{X}%
_{\alpha}\varphi(\mathbf{b}_{\alpha})\nonumber\\
&  =\delta_{\alpha}\text{.} \tag{17}%
\end{align}
We note that the computation of $||\mathbf{X}_{\alpha}||_{\infty\rightarrow1}$
is NP hard, see Rohn (2000)[$\mathbf{26}$]. The R package TaxicabCA by Allard
and Choulakian (2019)[$\mathbf{26}$]\ available online, uses three approaches
for its computation.

Besides equations (16) or (17), the parameters $(\delta_{\alpha}%
,\mathbf{a}_{\alpha},\mathbf{b}_{\alpha})$ for $\alpha=1,...,k$ satisfy:

a) The dispersion values are ordered; that is, $\delta_{1}\geq\delta_{2}%
\geq...\geq\delta_{k}>0.$

b) The sets $\mathbf{A}=[\mathbf{a}_{1}\ \mathbf{a}_{2}\ \mathbf{...}%
\ \mathbf{a}_{k}]\in R^{I\times k}$ and $\mathbf{B}=[\mathbf{b}_{1}%
\ \mathbf{b}_{2}\ \mathbf{...}\ \mathbf{b}_{k}]\in R^{J\times k}$ are $l_{p}$
- conjugate.

c) The transition formulas are for $\alpha=1,...,k$
\begin{align}
\mathbf{a}_{\alpha}  &  =\mathbf{X_{\alpha}\ }\varphi(\mathbf{b}_{\alpha
})\tag{18}\\
\text{\ }\mathbf{b}_{\alpha}  &  =\mathbf{X_{\alpha}^{\prime}\ }%
\varphi(\mathbf{a}_{\alpha})\nonumber\\
\delta_{\alpha}  &  =\varphi(\mathbf{a}_{\alpha})^{\prime}\mathbf{X_{\alpha
}\ }\varphi(\mathbf{b}_{\alpha}).\nonumber
\end{align}

Equations in (18), see Choulakian (2016)[$\mathbf{24}$], have simple
interpretations based on precedent sections: projection operator onto the
vector $\mathbf{a}_{\alpha}\in l_{p}^{I}$ is
\begin{align*}
\mathbf{Q}_{\mathbf{a}_{\alpha}}  &  \mathbf{=}\frac{\mathbf{a}_{\alpha}%
}{||\mathbf{a}_{\alpha}||_{p}}\varphi(\mathbf{a}_{\alpha})^{\prime}\\
&  =\frac{\mathbf{a}_{\alpha}}{\delta_{\alpha}}\varphi(\mathbf{a}_{\alpha
})^{\prime}\text{ by (16) and (17);}%
\end{align*}
and the projection of the $j$th column of $\mathbf{X_{\alpha}}$,
$\mathbf{X_{\alpha}(,}j\mathbf{)\in}l_{p}^{I}$ onto $\mathbf{a}_{\alpha}\in
l_{p}^{I}$ is%
\begin{align}
\mathbf{Q}_{\mathbf{a}_{\alpha}}\mathbf{X_{\alpha}(,}j\mathbf{)\ }  &
\mathbf{=\frac{\mathbf{a}_{\alpha}}{||\mathbf{a}_{\alpha}||_{p}}%
\varphi(\mathbf{a}_{\alpha})^{\prime}X_{\alpha}(,}j\mathbf{)}\tag{19}\\
&  =\mathbf{\ }b_{\alpha}(j)\mathbf{a}_{\alpha}/\delta_{\alpha}\text{ by
(18);}\nonumber
\end{align}
so the coordinate of the projection of the column $\mathbf{X_{\alpha}%
(,}j\mathbf{)}$ on the normed vector $\mathbf{a}_{\alpha}/\delta_{\alpha}\in
l_{p}^{I}$ is $b_{\alpha}(j)$. Thus based on these facts and designating the
$i$th row of $\mathbf{X}$ by $\mathbf{X(}i,\mathbf{),}$ we get the following
two results:

For $p=2$, based on the fact that the sets $\mathbf{A}=[\mathbf{a}%
_{1}\ \mathbf{a}_{2}\ \mathbf{...}\ \mathbf{a}_{k}]\in R^{I\times k}$ and
$\mathbf{B}=[\mathbf{b}_{1}\ \mathbf{b}_{2}\ \mathbf{...}\ \mathbf{b}_{k}]\in
R^{J\times k}$ are $l_{2}$ - conjugate (orthogonal), we have:%

\[
||\mathbf{X(,}j\mathbf{)\mathbf{||}}_{2}^{2}\ \mathbf{=||}\sum_{\alpha=1}%
^{k}b_{\alpha}(j)\mathbf{a}_{\alpha}/\delta_{\alpha}\mathbf{||}_{2}^{2}%
=\sum_{\alpha=1}^{k}(b_{\alpha}(j))^{2}.
\]

\[
||\mathbf{X(}i\mathbf{,)||}_{2}^{2}=\mathbf{||}\sum_{\alpha=1}^{k}a_{\alpha
}(i)\mathbf{b}_{\alpha}/\delta_{\alpha}\mathbf{||}_{2}^{2}=\sum_{\alpha=1}%
^{k}(a_{\alpha}(i))^{2}.
\]

\begin{equation}
||\mathbf{X\mathbf{||}}_{2}^{2}=\sum_{i,j}x_{ij}^{2}=\sum_{\alpha=1}^{k}%
\delta_{\alpha}^{2}. \tag{20}%
\end{equation}

For $p=1,$ based on the fact that the sets $\mathbf{A}=[\mathbf{a}%
_{1}\ \mathbf{a}_{2}\ \mathbf{...}\ \mathbf{a}_{k}]\in R^{I\times k}$ and
$\mathbf{B}=[\mathbf{b}_{1}\ \mathbf{b}_{2}\ \mathbf{...}\ \mathbf{b}_{k}]\in
R^{J\times k}$ are $l_{1}$ - conjugate, we have by applying the triangle
inequality of the $l_{1}$ - norm:%

\[
||\mathbf{X(,}j\mathbf{)||}_{1}=||\sum_{\alpha=1}^{k}b_{\alpha}(j)\mathbf{a}%
_{\alpha}/\delta_{\alpha}\mathbf{||}_{1}\leq\sum_{\alpha=1}^{k}|b_{\alpha
}(j)|.
\]

\[
||\mathbf{X(}i\mathbf{,)||}_{1}=||\sum_{\alpha=1}^{k}a_{\alpha}(i)\mathbf{b}%
_{\alpha}/\delta_{\alpha}\mathbf{||}_{1}\leq\sum_{\alpha=1}^{k}|a_{\alpha
}(i)|.
\]

\begin{equation}
||\mathbf{X\mathbf{||}}_{1}=\sum_{i,j}|x_{ij}|\ <\sum_{\alpha=1}^{k}%
\delta_{\alpha}\text{\ \ for\ \ }k>1\text{ \ by \ Corollary 1}. \tag{21}%
\end{equation}

\subsection{SVD and l$_{1}\min-$ SVD \ \ }

Let us consider the problem of
\begin{equation}
min\text{ \ }||\mathbf{X}_{\alpha+1}||_{p}^{p}=||\mathbf{X}_{\alpha
}-\mathbf{a}_{\alpha}\mathbf{b}_{\alpha}^{\prime}/\delta_{\alpha}||_{p}^{p}.
\tag{22}%
\end{equation}

For $p=2,$ $(22)$ corresponds to Pearson's (1901)[$\mathbf{20}$] approach of
PCA, that is, calculating SVD of \textbf{X} as in the precedent subsection.

For $p=1,$ to our knowledge, $(22)$ was first studied by Hawkins et al.
(2001)[$\mathbf{7}$] as a \textit{"Robust SVD"}, where they applied
\textit{"Alternating L1 regression"}, that is, alternating weighted median
method, see equation (12); see \ also among others, Ke and Kwak
(2005)[$\mathbf{28}$], Brooks and Dula (2019) [$\mathbf{29}$] and Song et al.
(2020)[$\mathbf{30}$]. In the last three cited papers, linear programming is
used to estimate the parameters in (22).

Hawkins et al. (2001)[$\mathbf{7}$] state that:

a) The dispersion values are no more ordered; that is, $\delta_{1}\geq
\delta_{2}\geq...\geq\delta_{k}>0$ is no more satisfied.

b) The relations among the columns of each matrix $\mathbf{A}=[\mathbf{a}%
_{1}\ \mathbf{a}_{2}\ \mathbf{...}\ \mathbf{a}_{k}]\in R^{I\times k}$ and
$\mathbf{B}=[\mathbf{b}_{1}\ \mathbf{b}_{2}\ \mathbf{...}\ \mathbf{b}_{k}]\in
R^{J\times k}$ are unknown. That is, matrices $\mathbf{A}$ and $\mathbf{B}$
are Grassman manifolds: $rank(\mathbf{A})=rank(\mathbf{B})=k$. While in the
ordinary SVD, matrices $\mathbf{A}$ and $\mathbf{B}$ are orthogonal,
equivalent to Stiefel manifolds $\mathbf{A}^{\prime}\mathbf{A}=\mathbf{B}%
^{\prime}\mathbf{B}=\mathbf{I}_{k}.$

c) The transition formulas are for $\alpha=1,...,k$
\begin{align}
a_{\alpha}(i)  &  =\arg\min_{\theta}||\mathbf{X_{\alpha}(}i\mathbf{,)-\ }%
\theta\mathbf{b}_{\alpha}||_{1}\tag{23}\\
b_{\alpha}(j)  &  =\arg\min_{\theta}||\mathbf{X_{\alpha}(,}j\mathbf{)-\ }%
\theta\mathbf{a}_{\alpha}||_{1}\nonumber\\
\delta_{\alpha}^{2}  &  =||\mathbf{a}_{\alpha}||_{1}||\mathbf{b}_{\alpha
}||_{1}.\nonumber
\end{align}

\subsection{Example 4}

Benz\'{e}cri (1973a,b) [$\mathbf{31a,31b}$], father of the French school of
data analysis, developed correspondence analysis (CA) based on weighted
Euclidean geometry. Benz\'{e}cri (1977, p 13) [$\mathbf{32}$] stated without
providing any argument: CA based on Euclidean geometry and on$\ l_{1}\min-$
SVD method are \textit{"qualitatively similar"}. The aim of this simple
example is to show that Benz\'{e}cri's statement seems to be true for
extremely sparse datasets.

Consider the diagonal matrix $\mathbf{X=Diag(}9\mathbf{\ \ \ \ }1).$ Then

a) SVD and $l_{1}\min-$ SVD of $\mathbf{X}$ are the same: $\mathbf{X}%
=9\mathbf{e}_{1}\mathbf{e}_{1}^{\prime}+\mathbf{e}_{2}\mathbf{e}_{2}^{\prime
},$ where $\mathbf{e}_{1}^{\prime}=(1\ \ \ 0)$ and $\mathbf{e}_{2}^{\prime
}=(0\ \ \ 1).$

\textit{Proof}: The solution of SVD is well-known. We provide the solution of
the $l_{1}\min-$ SVD by step-wise manner.%

\begin{align*}
min\ Loss(a,b)  &  =||\mathbf{X}-\mathbf{vv}^{\prime}||_{1}\text{ where
}\mathbf{v}^{\prime}=(a\ \ \ b)\text{ and }3\geq a\geq0,\text{ \ }1\geq
b\geq0\\
&  =|9-a^{2}|+2|ab|+|1-b^{2}|\ =9-a^{2}+2ab+1-b^{2}\\
&  =10-(a-b)^{2}\text{ \ and }3\geq a\geq0,\text{ \ }1\geq b\geq0\\
&  =10-3^{2}=1=Loss(a=3,b=0).
\end{align*}
We note that $Loss(a,b)=Loss(-a,-b);$ and also easily show that
$min\ Loss(a,b)\geq1$ for ($3\leq a,$ \ $1\leq b)$ and ($3\leq a,$ \ $1\geq
b\geq0)$ and ($3\geq a\geq0,$ \ $1\leq b).$ Thus $\mathbf{vv}^{\prime
}=9\mathbf{e}_{1}\mathbf{e}_{1}^{\prime}$ in step 1 and $\mathbf{vv}^{\prime
}=\mathbf{e}_{2}\mathbf{e}_{2}^{\prime}$ in step 2.

b) TSVD of $\mathbf{X}$\textbf{ }is computed in two steps: In step 1, given
that the two row vectors are located in the nonnegative quadrant, whose
representative is the corner vector $\mathbf{c}_{1}^{\prime}=(1\ \ 1),\ $so
the Taxicab projections of two row vectors onto $\mathbf{c}_{1}$ are the
coordinates of $\mathbf{u}_{1}=\mathbf{X(}_{1}^{1})=(_{1}^{9})$ and the
residual matrix is
\begin{align*}
\mathbf{X}_{2}  &  \mathbf{=}\left[
\begin{array}
[c]{cc}%
9 & 0\\
0 & 1
\end{array}
\right]  -\mathbf{u}_{1}\mathbf{u}_{1}^{\prime}/10\\
&  =\left[
\begin{array}
[c]{cc}%
0.9 & -0.9\\
-0.9 & 0.9
\end{array}
\right]  ;
\end{align*}
so $\mathbf{u}_{2}=\mathbf{X}_{2}\mathbf{(}_{-1}^{1})=\mathbf{(}_{-1.8}%
^{1.8})$ or $\mathbf{u}_{2}=\mathbf{X}_{2}\mathbf{(}_{1}^{-1})=\mathbf{(}%
_{1.8}^{-1.8})$; thus TSVD of $\mathbf{X}=\mathbf{u}_{1}\mathbf{u}_{1}%
^{\prime}/10+\mathbf{u}_{2}\mathbf{u}_{2}^{\prime}/3.6.$

\section{Geometric interpretation \textbf{of} the robustness \textbf{of} TSVD
}

\bigskip

\subsection{Curse of dimensionality}

The curse of dimensionality introduced by Bellman (1961)\ [$\mathbf{33}$] is a
central issue in modern data analysis, particularly in the study and
visualization of high-dimensional datasets. Let $\mathbf{X}=(x_{ij})$ denote
data matrix of size $n\times d$. In the literature, a dataset is typically
considered high-dimensional in two main situations:

Case a) The number of variables exceeds the sample size $d>n$, where $d$
denotes the number of continuous variables and $n$ the number of observations;
see, for example, Pires and Branco (2019) [$\mathbf{34}$].

Case b) The data are sparse or extremely sparse, a situation that is
increasingly common in contemporary applications.

In the remainder of this section, \textit{we assume that the data are
uniformly distributed over both the hypersphere and the hypercube}. The
effectiveness of TSVD in analyzing and visualizing high-dimensional data is
rooted in a well-known geometric property of such spaces: as the
dimensionality $d$ increases, the volume of the hypercube becomes increasingly
concentrated near its corners $(\pm1)^{d}$, while the volume of the largest
inscribed hypersphere tends to zero.

This phenomenon has important implications for the computation and behavior of
principal axes, which we explore in detail. A thorough discussion of these
geometric aspects can be found, among others, in Scott (2015, Section 1.5)
[$\mathbf{15}$] and Zaki and Meira (2014, Chapter 6) [$\mathbf{16}$], who cite
Kendall (1961) [$\mathbf{17}$] as a primary reference.

\subsection{\textbf{Volumes of Hypersphere and Hypercube }}

Let $a=max_{i,j}|x_{ij}|$ for $i=1,...,n$ and $j=1,...,d.$ The $n$ data points
are located both

a) within the $d$-dimensional hypercube%

\begin{align*}
H_{\infty}^{d}(a)  &  =\{\mathbf{y}:||\mathbf{y||}_{\infty}=\max
\ |y_{i}|\ =a\text{ for }i=1,...,d\}\\
&  =[-a,a]^{d}%
\end{align*}
whose volume is%

\[
V_{\infty}^{d}(a)=(2a)^{d}.
\]

Thus%
\begin{align*}
V_{\infty}^{d}(a  &  >1/2)\rightarrow\infty\text{ as }d\rightarrow\infty\\
V_{\infty}^{d}(a  &  =1/2)\rightarrow1\text{ as }d\rightarrow\infty\\
V_{\infty}^{d}(a  &  <1/2.)\rightarrow0\text{ as }d\rightarrow\infty.
\end{align*}

b) within the $d$-dimensional Euclidean hypersphere%

\[
H_{2}^{d}(a)=\{\mathbf{y}:||\mathbf{y}||_{2}^{2}=\sum_{j=1}^{d}y_{j}^{2}%
=a^{2}\},
\]
whose volume is%

\[
V_{2}^{d}(a)=\frac{a^{d}\ \pi^{d/2}}{\ \Gamma(\frac{d}{2}+1)},
\]
$\ $where $\Gamma(\frac{d}{2}+1)$ is the Gamma function. Stirling's
approximation of the Gamma function is%

\[
\Gamma(x+1)=x!\simeq x^{x}e^{-x}(2\ \pi\ x)^{1/2}.
\]

It is evident that the hypersphere $H_{2}^{d}(a)$ is inscribed within the
hypercube $H_{\infty}^{d}(a)$; that is, $V_{2}^{d}(a)<V_{\infty}^{d}(a)$; so
the ratio of the two volumes is given by the fraction%

\begin{align*}
f_{d}  &  =\frac{V_{2}^{d}(a)}{V_{\infty}^{d}(a)}=\frac{a^{d}\ \pi^{d/2}%
}{(2a)^{d}\ \ \Gamma(\frac{d}{2}+1)}\\
&  =\frac{\ \pi^{d/2}}{2^{d}\ \ \Gamma(\frac{d}{2}+1)}=\frac{(\pi/4)^{d/2}%
}{\Gamma(\frac{d}{2}+1)}\\
&  <\frac{1}{\Gamma(\frac{d}{2}+1)};
\end{align*}
thus
\[
f_{d}\rightarrow0\text{ as }d\rightarrow\infty;
\]
from which we deduce%
\[
V_{2}^{d}(a)\rightarrow0\text{ as }d\rightarrow\infty\text{ for }a>0.
\]

Scott (2015, page 30, Table 1.1)[$\mathbf{15}$] provides two values
$f_{d=3}=0.524$ and $f_{d=7}=0.037.$

\subsubsection{\textbf{Volume of a hypercube near its corners }}

It is evident that the center of the hypercube becomes progressively less
relevant as the dimension $d$ increases. In high dimensions, the volume of the
hypercube concentrates near its corners $(\pm1)^{d}$. For high-dimensional
datasets, particularly in the context of dimension reduction for
visualization, this geometric phenomenon helps explain the robustness of TSVD.
Indeed, under Grothendieck's norm (see equation (17)), a taxicab principal
axis corresponds to the direction defined by a corner $\mathbf{c}$ and its
opposite $\mathbf{-c}$, thereby capturing the dominant structure of the data.

\subsubsection{\textbf{Volume of a hypersphere in a thin shell }}

For high-dimensional datasets, Wegman (1990)[$\mathbf{35}$] demonstrated the
following: Consider two centered hyperspheres, one with radius $a$ and the
other with slightly smaller radius $a-\epsilon$. Then,%

\[
\frac{V_{2}^{d}(a)-V_{2}^{d}(a-\epsilon)}{V_{2}^{d}(a)}=\frac{a^{_{^{^{d}}}%
}-(a-\epsilon)^{d}}{a^{_{^{^{d}}}}}=1-\frac{(a-\epsilon)^{d}}{a^{_{^{^{d}}}}%
}\rightarrow1\text{ \ as}\ d\rightarrow\infty
\]

Thus, for high-dimensional datasets virtually all of the volume of a
hypersphere is concentrated in a thin shell close to its surface, which is
only a $(d-1)$-dimensional manifold.

The next result is very common in applications.

\subsubsection{\textbf{The angle }$\theta_{d}$\textbf{ between a corner vector
}$\mathbf{c}$\textbf{ and a particular vector }}

Let $\mathbf{c=(1}_{k}\mathbf{,\pm1}_{d-k})^{\prime}$ for $1\leq k\ll d$ be a
hypercube corner vector, and $\mathbf{a=(1}_{k}\mathbf{,0}_{d-k})^{\prime}$ a
combination of $k$ standard unit vectors such as $\mathbf{e}_{1}%
\mathbf{=(}1,\mathbf{0}_{d-1})^{\prime}$, $\mathbf{e}_{j}\mathbf{=(0}%
_{j-1},1,\mathbf{0}_{d-j})^{\prime}$ for $j=2,...,k.$ Then $\mathbf{a}$
represents both the orthogonal (Euclidean) and $l_{1\min}$ projection of
$\mathbf{c}$ onto $\mathbf{a}.$ Thus by the two Pythagorean Theorems
\begin{align*}
||\mathbf{c}||_{p}^{p}  &  =||\mathbf{a}||_{p}^{p}+||\mathbf{c-a}||_{p}%
^{p}\text{ \ \ for }p=1,2\\
d  &  =k+(d-k);
\end{align*}
so, in the Euclidean geometry $cos(\theta_{d})=||\mathbf{a}||_{2}%
/||\mathbf{c}||_{2}=\sqrt{k/d}\rightarrow0$ as $d\rightarrow\infty;$ similarly
in Taxicab geometry $cos(\theta_{d})=||\mathbf{a}||_{1}/||\mathbf{c}%
||_{1}=k/d\rightarrow0$ as $d\rightarrow\infty$; that is, for large $d$ the
angle $\theta_{d}\asymp90^{0}$ between a corner vector $\mathbf{c}$ and the
particular vector $\mathbf{a}$. This is paradoxical at first glance, although
it is a well-known phenomenon in the behavior of SVD in high-dimensional
datasets. Accordingly, we obtain the following two commonly observed
properties for high-dimensional sparse or extremely sparse datasets, where the
$n$ data points in $R^{d}$ are uniformly distributed:

a) The vector $\mathbf{a}$ represents an approximate principal axis in both
SVD and $l_{1\min}$-SVD, on which only a finite number of projected points
$n_{1}\geq$ $k$\ lie. Since $k\ll d$ and $d$ is large, it follows from the
result in subsection 7.2.1 the remaining $(n-n_{1})$ points are located near
the hypercube corners $\mathbf{c}$. Their orthogonal or $l_{1\min}$
projections onto the principal axis $\mathbf{a}$ lie close to the origin due
to the angle $\theta_{d}\asymp90^{0}.$ Example 4 is a simple case study of
this phenomenon.

(b) In TSVD, a principal axis coincides with an orthant corner vector
$\mathbf{a}=\mathbf{c}_{1}$. Therefore, for uniformly distributed data of size
$n$, the projections onto the principal axis $\mathbf{a}=\mathbf{c}_{1}$ are
themselves uniformly distributed.

\section{Conclusion}

We analyzed and compared three modern statistical approaches to dimensionality
reduction based on stepwise double projections, highlighting their Euclidean
and Taxicab geometric foundations, as well as their effectiveness and
practical relevance---particularly for high-dimensional and extremely sparse data.

In Taxicab geometry, a centered space of dimension $d$ consists of $2^{d}$
orthants, each represented by a unit corner vector $\mathbf{c}$ whose
coordinates take values of $\pm1$. Consequently, in TSVD, all principal axes
are unit corner vectors---an aspect that helps explain its robustness.

We conclude by emphasizing the value of jointly using SVD and TSVD. Comparing
the visual maps produced by both methods is strongly recommended. Choulakian
and Allard (2025)  [$\mathbf{9}$] explored the combined application of
marginally weighted SVD (CA) and its robust counterpart, marginally weighted
TSVD (TCA), for the analysis and visualization of an extremely sparse,
high-dimensional dataset. Within this framework, TSVD can be viewed either as
a robust alternative to classical SVD or as a complementary approach, with
both methods offering distinct yet mutually informative perspectives on the
underlying structure of the data.\bigskip

\textbf{Declarations}

\textit{Author Contribution}: Vartan Choulakian conceptualization and writing.

\textit{Funding}: Partial funding to the author is provided by the Natural
Sciences and Engineering Research of Canada (Grant no. RGPIN-2017-05092).

\textit{Data availability}: None.

\textit{Conflict of interest}: There is no conflict of interest.

\bigskip

\bigskip\textbf{References ordered}

\textbf{1} Agniel V (2021) Lp-Projections on Subspaces and Quotients of Banach
Spaces. \textit{Advances in Operator Theory}, 6,

https://doi.org/10.1007/s43036-021-00131-8

\textbf{2} Alber Y (1996) Metric and generalized projection operators in
Banach spaces: Properties and applications. In \textit{Theory and Applications
of Nonlinear Operators of Monotone and Accretive `Qpe}, (Edited by A.
Kartsatos), pp. 15-50, Marcel Dekker, New York

\textbf{3} Alber Y (2017) Generalized projections and equivalent
representations of James orthogonal decompositions in Banach spaces.
\textit{Communications on Applied Nonlinear Analysis}, 24 (2), 28-48

\textbf{4} Golland L (1990) Karl Menger and Taxicab Geometry.
\textit{Mathematics Magazine}, 63 (5): 326--327. doi:10.1080/0025570x.1990.11977548

\textbf{5} Eckart C, Young G (1936) The approximation of one matrix by another
of lower rank. \textit{Psychometrika}, 1, 211--218

\textbf{6} Choulakian V (2006) Taxicab correspondence analysis.
\textit{Psychometrika,} 71, 333-345

\textbf{7} Hawkins DM, Liu L, Young SS (2001) Robust singular value
decomposition. \textit{NISS Technical Report 122.}

\textit{ }Available at (www.niss.org/downloadabletechreports.html)

\textbf{8} Choulakian V (2017) Taxicab correspondence analysis of sparse
contingency tables. \textit{Italian Journal of Applied Statistics,} 29 (2-3), 153-179

\textbf{9} Choulakian V, Allard J (2025) Scale-invariant correspondence
analysis of compositional data. \textit{Analytics}, 4, 32. Available at

https://doi.org/10.3390/analytics4040032

\textbf{10} Gauthier SM, Choulakian V (2015) Taxicab correspondence analysis
of abundance data in archeology: three case studies revisited.
\textit{Archeologia e Calcolatori}, 26, 77-94

\textbf{11} Choulakian V, Tibeiro JD (2013a) Graph partitioning by
correspondence analysis and taxicab correspondence analysis. \textit{Journal
of} \textit{Classification} 30:397-427

\textbf{12} Choulakian V, Kasparian S, Miyake M, Akama H, Makoshi N, Nakagawa
M (2006) A Statistical Analysis of the Synoptic Gospels. \textit{Journ\'{e}es}
\textit{Internationales d'Analyse Statistique des Donn\'{e}es Textuelles} 8: 281-88

\textbf{13} Choulakian, V., Allard, J. and Simonetti, B. (2013b) Multiple
Taxicab correspondence analysis of a survey related to health services.
\textit{Journal of Data Science}, 11, 205-229.

\textbf{14} Choulakian V, Allard J, Smail M (2023) Taxicab correspondence
analysis and taxicab logratio analysis: a comparison on contingency tables and
compositional data. \textit{Austrian Journal of Statistics}, 52(3), 39-70

\textbf{15} Scott DW (2015) \textit{Multivariate density estimation: theory,
practice, and visualization}. Second edition, John Wiley \& Sons, Inc

\textbf{16} Zaki MJ, Meira WJ (2014) Data Mining and Analysis: Fundamental
Concepts and Algorithms. Cambridge University Press

\textbf{17} Kendall MG (1961) \textit{A Course in the Geometry of n
Dimensions}. Hafner Publishing Company, New York, USA

\textbf{18} Choulakian V (2016) Matrix factorizations based on induced norms.
\textit{Statistics, Optimization and Information Computing}, 4, 1-14

\textbf{19} Gardner M (1997) \textit{The Last Recreations}. Springer-Verlag
New York, Inc

\textbf{20} Pearson K (1901) On lines and planes of closest fit to systems of
points in space. \textit{Philosophical Magazine,} 2:559-572

\textbf{21} Hotelling H (1933) Analysis of a complex of statistical variables
into principal components. \textit{Journal of Educational Psychology}, 24, 417-441

\textbf{22} Jolliffe I (2002) \textit{Principal component analysis}. Springer
New York, NY

\textbf{23} Choulakian V (2003) The optimality of the centroid method.
\textit{Psychometrika}, 68(3), 473-475

\textbf{24} Burt C (1917) \textit{The Distribution and Relations of
Educational Abilities}. P.S. King \& Son, London

\textbf{25 }Thurstone LL (1931) Multiple factor analysis.
\textit{Psychological Rev}. 38, 406--427

\textbf{26} Rohn J (2000) Computing the norm $\parallel\mathbf{A}%
\parallel_{\infty,1}$ is NP-hard. \textit{Linear and Multilinear Algebra}, 47
(3), 195-204

\textbf{27} Allard J, Choulakian V (2019) \textit{Package TaxicabCA in R}

\textbf{28} Ke, Kanade T (2005) Robust L1 norm factorization in the presence
of outliers and missing data by alternative convex programming. In
\textit{Proceedings of the IEEE Conference on Computer Vision and Pattern
Recognition}

\textbf{29} Brooks JP, Dula JH (2019) Approximating L1-norm best-fit lines.
\textit{Optimization Online}

\textbf{30} Song Z, Woodruff DP, Zhong P (2020) Low Rank Approximation with
Entrywise $l_{1}-$Norm Error. Available at:

\textit{http://arxiv.org/abs/1611.00898}

\textbf{31a} Benz\'{e}cri JP (1973a) \textit{L'Analyse des Donn\'{e}es: Vol.
1: La Taxinomie}. Paris: Dunod.

\textbf{31b} Benz\'{e}cri JP (1973b) \textit{L'Analyse des Donn\'{e}es: Vol.
2:} \textit{L'Analyse des Correspondances}. Paris: Dunod.

\textbf{32} Benz\'{e}cri JP (1977) Histoire et Pr\'{e}histoire de l'Analyse
des Donn\'{e}es: V: L'Analyse des correspondances. \textit{Les Cahiers de
l'Analyse des} \textit{Donn\'{e}es}, II(1), 9-53

\textbf{33} Bellman RE (1961) Adaptive Control Processes. \textit{Princeton
University Press}, Princeton, NJ

\textbf{34} Pires AC, Branco JA (2019) Branco. High dimensionality: The latest
challenge to data analysis. Available at:

arXiv:1905.02330

\textbf{35 }Wegman EJ (1990) Hyperdimensional Data Analysis Using Parallel
Coordinates. \textit{Journal American Statististical Association,} 85 664--675

\end{document}